\documentclass[11pt]{article}
\usepackage[cp1251]{inputenc}
\usepackage[english,russian]{babel}
\usepackage{amsfonts,amssymb,eucal,amsmath}
\usepackage{epic,graphicx,bezier,curvesls}
\usepackage{latexsym}
\usepackage{theorem}

\theoremstyle{plain} %[section]

\begin{document}

\begin{center}
{\large\textbf{FINITE SYSTEMS OF EQUATIONS AND IMPLICIT FUNCTIONS} \\[3mm]%
\textbf{P.P. Zabreiko, A.V. Krivko-Krasko}} \\

\vspace{3mm}
\end{center}

%%%%%%%%% Introduction %%%%%%%%%

The classical R\"uckert--Lefschetz scheme of analysis of implicit functions (defined by finite
systems of $n$ analytical equations with $n$ unknowns) is studied from the point of view of
calculations with finite number coefficients in Taylor expansions for left hand parts of
corresponding equations. It is proved that this scheme is not applicable in the general case. It is
offered some modifications allowing to lead the analysis of implicit functions to the calculation
with finite number of coefficients.

%%%%%%%%% Introduction %%%%%%%%%

{\bf Introduction.} Let us consider a finite system of the equations
 \begin{equation}\label{scsyst}
 \left\{ \begin{array}{c}f_1(\lambda,x_1,\ldots,x_m) = 0, \\ \ldots \ldots \ldots \ldots \ldots \ldots \ldots \\
 f_n(\lambda,x_1,\ldots,x_m) = 0,\end{array} \right.
 \end{equation}
where the parameter $\lambda$ and the unknowns $x_1,\ldots,x_m$ are real or complex numbers and
$f_j(\lambda,x_1,$ $\ldots,x_m)$ ($j = 1,\ldots,n$) are real or complex valued functions. System
(\ref{scsyst}) can be written as
 \begin{equation}\label{vecsyst}
 f(\lambda,x) = 0,
 \end{equation}
where $f(\cdot,\cdot)$ is a map of ${\Bbb R}\times {\Bbb R}^m$ to ${\Bbb R}^m$ or ${\Bbb C}\times
{\Bbb C}^m$ to ${\Bbb C}^m$.

Suppose that
 $$
 f(\lambda_0,x_0) = 0.
 $$
In a number of problems of analysis (in particular, in differential and integral equations,
optimization methods and etc) the following question arises: when does System (\ref{vecsyst})
define in a neighborhood of the point $(\lambda_0,x_0)$ (or in some part of this neighborhood) one
or several functions $x(\lambda)$ that are continuous at the point $\lambda_0$ and such that
$x(\lambda_0) = x_0$? This functions are often called {\it implicit functions} or {\it small
solutions} of System (\ref{scsyst}).

The classical theorem about implicit functions is well-known \cite {KVZRS, FIHT}: {\it if $m = n$,
$f (\lambda,x_0)$ is a continuous function at the point $\lambda_0$, $f'_x(\lambda,x)$ is a
continuous function at the point $(\lambda_0,x_0)$ and $f'_x(\lambda_0,x_0)^{-1}$ exists, then
System {\rm (\ref{vecsyst})} has a unique solution $x = x^*(\lambda)$ in a small neighborhood of
the point $(\lambda_0,x_0)$}. This case is called {\it nondegenerated}. If $f'_x (\lambda_0,x_0)$
is an irreversible matrix, then the corresponding case is called {\it degenerated}.

The analysis of degenerated cases is a difficult problem. The basic results concern the case when
$f_j(\lambda,x_1,\ldots,x_m)$ ($j = 1,\ldots,n$) are analytical functions in a neighborhood of
$(\lambda_0,x_0)$ (see, for example, \cite {VT, E, KVZRS}); some of these results are extended to
the case when the functions $f_j(\lambda,x_1,\ldots,x_m)$ are smooth enough.

Depending on which of the cases $m = n$, $m > n$, $m < n$ takes place, it is said that System
(\ref{scsyst}) is {\it determined}, {\it underdetermined} and {\it overdetermined}. It seems that
determined systems should define a finite number of solutions $x(\lambda)$, underdetermined ones
should define infinite number of such solutions and overdetermined ones should not define any
solution in general. However the distinction between these three types of systems is conditional.
So if we add one or several equation so that the number of equations became the same as the number
of unknowns then the underdetermined system becomes determined. The overdetermined system also can
be considered as determined system if the left hand parts of its equations depend on also $n - m$
additional unknowns $x_{m+1},\ldots,x_n$.

In the article (if it is not stipulated the opposite) we consider the case when the parameter
$\lambda$ and the unknowns $x_1,\ldots,x_m$ take complex values. There are situations when the
solutions in question <<branch>> at the point $\lambda_0$. To avoid consideration of
multiple-valued functions in such cases it is natural to consider the implicit functions defined by
System (\ref{scsyst}) in the neighborhood of the point $\lambda_0$ with a cross-cut. The case when
the parameter $\lambda $ and unknowns $x_1,\ldots,x_m$ take real values will be considered in
detail in the second part of this article.

Assume that the functions $f_j(\lambda,x_1,\ldots,x_m)$ ($j = 1,\ldots,n$) are analytical. Then the
zero set ${\frak N} = \{(\lambda, x) \in {\Bbb C}^{m+1}: \ \lambda \ne \lambda_0, \ f(\lambda,x) =
0 \}$ of the left hand parts of System (\ref{scsyst}) in the neighborhood of the point
$(\lambda_0,x_0)$ can be presented in the form ${\frak N} = {\frak N}_0 \cup {\frak N}_1 \ldots
\cup {\frak N}_m$. Here the set ${\frak N}_0$ is empty or consists of a finite number of graphs of
solutions $x = \phi(\lambda)$ where $\phi(\lambda)$ are some analytical functions of the parameter
$(\lambda - \lambda_0)^\frac1r$ ($r$ is a natural number). Further, each of the sets ${\frak N}_j$
($j = 1,\ldots,m - 1$) is empty or consists of a finite number of the <<surfaces>> that are the
graphs of functions of type $x = \phi(\lambda;\xi_1,\ldots,\xi_j)$ where $\xi_s$ ($s = 1,\ldots,j$)
are free parameters (fig. 1).
 \begin {figure} [h]
 \centering
 \includegraphics [scale=0.7] {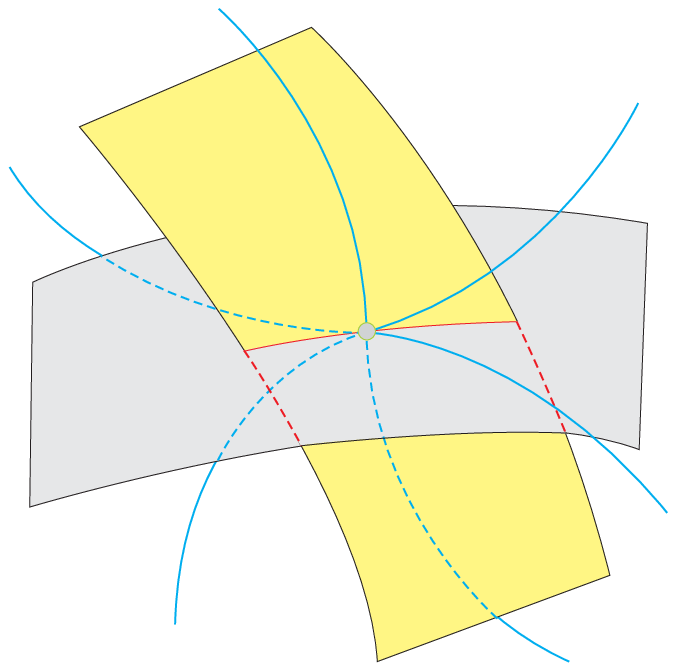} \\
 Fig. 1. The zero set ${\frak N}$
 \end {figure}
Moreover, the functions $\phi(\lambda;\xi_1,\ldots,\xi_j)$ ($j = 1,\ldots, m - 1$) have the
following property: if we replace in these functions the parameters $\xi_1,\ldots,\xi_j$ ($s =
1,\ldots,j$) by some analytical functions $\xi_s(\lambda)$ ($s = 1,\ldots,j$) depending on
parameters $(\lambda - \lambda_0)^\frac1p$ ($p$ is a natural) then the superpositions
$\phi(\lambda;\xi_1(\lambda),\ldots,\xi_j(\lambda))$ also will be analytical functions of the
parameter $(\lambda - \lambda_0)^\frac1q$ ($q$ is also a natural, and $p$ is a divisor of $q$). At
last, the set ${\frak N}_m$ is not empty only when System (\ref{scsyst}) is trivial, i. e. when its
left hand parts are identically equal to zero (in this case any continuous at the point $\lambda_0$
function $x(\lambda)$ such that $x(\lambda_0) = x_0$ satisfies System (\ref{scsyst})).

As a result we can give a description of the general structure of implicit functions defined by
System (\ref{scsyst}) by means of the objects mentioned above. Obviously the function $x =
x(\lambda)$ defined in the neighborhood of the point $\lambda_0$ is implicit if and only if its
graph lies in the set ${\frak N}$. In particular, if the set ${\frak N}_0$ is not empty then it
defines a finite set $\Phi_0$ of implicit functions. Each of these functions is an analytical
function of the parameter $(\lambda - \lambda_0)^\frac1r $ ($r$ is a natural number). Further, the
sets $\Phi_j$ ($j = 1,\ldots,m - 1$) of implicit functions, that are analytical functions of the
parameter $(\lambda - \lambda_0)^\frac1r$ ($r$ is a natural number), and whose graphs lie on the
<<surfaces>> <<consistuting >> the set ${\frak N}_j$ ($j = 1,\ldots,m - 1$), are infinite provided
that they are nonempty. Remark also the following: if the set ${\frak N}_j$ ($j = 1,\ldots,m-1$) is
not empty then there are others (continuous in a neighborhood $\lambda_0$!) implicit functions
which are not analytical of the parameter $(\lambda - \lambda_0)^\frac1r$ ($r$ is a natural
number). However such functions can be excluded from the consideration as soon as the graphs of the
analytical implicit functions of the set $\Phi_j$ fill the surfaces of ${\frak N}_j$ in the
neighborhoods of the point $(\lambda_0,x_0)$.

Let us denote by $\Phi$ the set of all analytical implicit functions $x = x(\lambda)$ of the
parameters $(\lambda - \lambda_0)^\frac1r$ ($r$ is a natural number) taking the value $x_0$ at
$\lambda = \lambda_0$. Obviously $\Phi = \Phi_0 \cup \Phi_1 \cup \ldots \cup \Phi_m$.

As far as we know, the above formulated statements were proved in the first half of \mbox{XX-th}
century by W. R\"uckert \cite{R}. The more modern statement of these results can be found in the
monographs \cite{BM,GR} (see also, \cite {KVZRS}). The corresponding argument was based on
Kronecker elemination theory for systems of algebraic equations and on Weierstrass preparation
theorem for analytical functions of the complex variable. In the monograph by S. Lefschetz \cite{L}
(he investigated special systems of type (\ref {scsyst}) which arose in the problem about periodic
solutions of ordinary differential equations) the more elementary statement of the results about
the structure of the set of the implicit functions was given.

Though S. Lefschetz's argument was not constructive, it is laid down in the basis of the general
constructions of M. M. Vainberg and V. A. Trenogin. In the monograph \cite{VT} they stated that
their scheme allows them to give the complete description of the sets $\Phi_0,\Phi_1,\ldots,\Phi_m$
and, moreover, to define the first coefficients in the expansions of solutions from the set
$\Phi_0$ into series along the parameter $\lambda$ or along its fractional degrees. In the
monograph \cite{KVZRS} it was noticed that it is not true. In this monograph it was shown that the
first coefficients of the expansions of solutions of the general system (\ref{scsyst}) into series
can be defined probably only for the so-called simple solutions (the solution $\phi(\lambda)$ of
System (\ref{scsyst}) is {\it simple} if $m = n$ and for the values $\lambda$ that are close to
$\lambda_0$ and distinct from $\lambda_0$ the Jacobian $\det f'_x (\lambda,\phi(\lambda))$ is non
zero). Moreover, in this monograph it was shown that the scheme of M.~ M.~ Vainberg and V. A.
Trenogin does not allow (if we use in the calculations only a finite number of coefficients of the
expansions into series of the left hand part of System (\ref{scsyst})) to define the number of the
implicit functions of the set $\Phi_0$ and the coefficients of the first members of the expansions
of these implicit functions even when $m = n = 3$ (and in essence when $m = n = 2$). In this
monograph it was also given a special example of the system (when $m = n = 3$) when some updating
of the Lefschetz scheme allow us to define the structure of the set $\Phi_0$.

\smallskip

In the next section the R\"uckert--Lefschetz scheme will be analyzed in detail. Besides in this
section we emphasize some moments which make the R\"uckert--Lefschetz scheme <<not constructive>>
and, moreover, the R\"uckert--Lefschetz scheme does not allow to define the structure of the set of
the implicit functions defined by System (\ref{scsyst}) even for the rough systems (in this article
System (\ref{scsyst}) is called {\it rough} if it has only a finite number of simple solutions (in
\cite {KVZRS} the term {\it rough systems} was used in a bit different sense). It is known that
System (\ref{scsyst}) is rough if and only if it possesses the following {\it stability property}:
for every big enough natural $N$ there exists a natural $\widetilde {N}$ such that if we change the
members in the left hand parts of System (\ref{scsyst}) whose orders are higher than
$\widetilde{N}$ then the number of the solutions of System (\ref{scsyst}) does not change and,
moreover, the first $N$ members of the expansions of these solutions into series also do not
change. In the fourth section a certain modified scheme of the research of System (\ref{scsyst}) is
offered; the basic idea of this modification is due to the above mentioned example from
\cite{KVZRS}.

Let us notice that the R\"uckert--Lefschetz scheme is not unique. The various statements about the
structure of the implicit functions defined by System (\ref{scsyst}) have been received by V. V.
Pokornyi, P. P. Rybin, V. B. Melamed, A. E. Gel'man; the considerable part of the results of these
authors is summarized in the monograph \cite{KVZRS}. It is necessary to note separately the
monograph \cite{E} by N. P. Erugin because his work contain a number of theorems about implicit
functions which are based on the construction of jets (the sums of first members in expansions into
power series) of the expansions of these functions into series.

The proof of the R\"uckert--Lefschetz scheme use only the elementary means of algebra and the
theory of functions of complex variables (i. e., such classical concepts as the resultant, the
greatest common divisor of polynomials with coefficients from the factorial rings (i. e. the rings
with the unique factorization on primes), etc.). The abstract theory of the polynomial ideals is
not used. In this article we use the results on the theory of implicit functions which are
described in \cite {KVZRS}.

\medskip

{\bf The R\"uckert--Lefschetz scheme.} Below we assume that $\lambda_0 = 0$, $x_0 = 0$ and $f_j
(\lambda, x_1, \ldots, x_m)$ $\not\equiv 0$ ($j = 1,\ldots,n$). The latter assumption implies
$\Phi_m = \emptyset$.

We change the notation of functions $f_j(\lambda,x_1,\ldots, x_m)$ onto
$f^{(m)}_j(\lambda,x_1,\ldots,x_m)$ $(j = 1,\ldots,n)$ (in what follows, it is convenient to fix
the number of unknowns in designations). Since the functions $f^{(m)}_j(\lambda,x_1,\ldots,x_m)$,
$j = 1,\ldots,n$, are analytical, we present the functions $f^{(m)}_j(\lambda,x_1,\ldots,x_m)$ $(j
= 1,\ldots,n)$ in the form of converging series in some neighborhood of zero
 \begin{equation}\label{row}
 f^{(m)}_j(\lambda,x_1,\ldots,x_m) = \sum _ {k_0 + k_1 + \ldots + k_m = 1} ^ \infty a _ {k_0, k_1, \ldots, k_m}
 \lambda ^ {k_0} x_1 ^ {k_1} \ldots x_m ^ {k_m} \qquad (j = 1, \ldots, n).
 \end{equation}

We divide each equation of System (\ref{row}) by the highest possible degree of $\lambda$ and so,
without the loss of generality, we can assume
 $$f_j^{(m)}(0,x_1,\ldots,x_m) \not\equiv 0 \quad (j = 1,\ldots,n).$$
In addition we make a linear substitution of the unknowns $x_1,\ldots,x_m$ so that the functions
 \begin{equation}\label{regu}
 f_j^{(m)}(0,0,\ldots,0,x_m) \qquad (j = 1,\ldots,n)
 \end{equation}
turn out to be nonzero.

As a result of the application of the Weierstrass preparation theorem \cite{GR} to each function
$f_j^{(m)}(\lambda,x_1,\ldots,x_m)$ ($j = 1,\ldots,n$) we receive the equalities
 \begin{equation}\label{veier}
 f_j^{(m)}(\lambda,x_1,\ldots,x_m) = \varepsilon^{(m)}_j(\lambda,x_1,\ldots,x_m) \cdot
 \widetilde{f}^{(m)}_j(\lambda,x_1,\ldots,x_m),
 \end{equation}
where $\varepsilon^{(m)}_j(\cdot)$ is an analytical function at zero, such that
$\varepsilon^{(m)}_j(0) \neq 0$; $\widetilde{f}^{(m)}_j(\lambda,x_1,\ldots,x_m)$ is a polynomial
with respect to the unknown $x_m$ whose coefficients are analytical at zero functions of the
parameter $\lambda $ and unknowns $x_1,\ldots,x_{m-1}$.

The equalities (\ref{veier}) imply that the search of implicit functions defined by System
(\ref{scsyst}) is equivalent to the analysis of the system of algebraic equations with respect to
the unknown $x_m$:
 \begin{equation}\label{newsys}
 \widetilde{f}^{(m)}_j(\lambda,x_1,\ldots,x_m) = 0 \qquad (j = 1,\ldots,n).
 \end{equation}
Notice that the superior coefficients of the polynomials of $x_m$ in the left hand parts of System
(\ref{newsys}) are equal to $1$ and all the other coefficients of these polynomials are analytical
functions of the parameter $\lambda$ and unknowns $x_1,\ldots,x_{m-1}$ turning into zero at zero.
This follows from our construction.

Let us remind (see, for example, \cite{KVZRS,L,VW}) that the rings of analytical at zero functions
of a finite number of variables are factorial (a ring is called factorial if it has an identity
element, has no divisors of zero, and its elements are (uniquely up to the order of multipliers)
displayed as the product of prime multipliers). In such rings the concept of the greatest common
divisor is defined and all the main statements of the divisibility theory are true. In particular,
the ring of polynomials with coefficients from a factorial ring itself is a factorial ring.

Let us denote by $d^{(m)}(\lambda,x_1,\ldots,x_m)$ the greatest common divisor of polynomials
$\widetilde{f}^{(m)}_j(\lambda,x_1,\ldots,$ $x_m)$, $j = 1,\ldots,n$. Then
 \begin{equation}\label{nodd}
 \begin {array} {c}
 \widetilde{f}^{(m)}_j(\lambda,x_1,\ldots,x_m) = d^{(m)}(\lambda,x_1,\ldots,x_m) \cdot
 \widetilde{\widetilde{f}}^{(m)}_j(\lambda,x_1,\ldots,x_m) \\ (j = 1,\ldots,n).
 \end{array}
 \end{equation}
Hence System (\ref{scsyst}) is equivalent to the collection consisting of one algebraic equation
with the unknown $x_m$
 \begin{equation}\label{nod}
 d^{(m)}(\lambda,x_1,\ldots,x_m) = 0,
 \end{equation}
and the system of the algebraic equations with the unknown $x_m$
 \begin{equation}\label {syx}
 \widetilde{\widetilde{f}}^{(m)}_j(\lambda,x_1,\ldots,x_m) = 0 \qquad (j = 1,\ldots,n).
 \end{equation}

If the degree $s$ of the greatest common divisor $d^{(m)}(\lambda,x_1,\ldots,x_m)$ is positive then
Equation (\ref{nod}), for any small enough $\lambda,x_1,\ldots,x_{m-1}$, has $s$ small solutions
$x_m$. These solutions can be presented as the equations $x_m = \phi(\lambda,x_1,\ldots,x_{m-1})$.
More precisely, each of such equations defines an element of the set $\Phi_{m-1}$ and the solutions
$x(\lambda) = (x_1(\lambda),\ldots,x_{m-1}(\lambda),x_m(\lambda))$ of System (\ref{scsyst}) where
$x_i(\lambda)$ ($i = 1,\ldots,m - 1$) are arbitrary analytical functions of $\lambda$ or of some
fractional degree $\lambda^\frac1r$ of $\lambda$ turning into zero at zero and the component
$x_m(\lambda)$ is defined by the equation
 $$x_m(\lambda) = \phi(\lambda;x_1(\lambda),\ldots,x_{m-1}(\lambda))$$
where $\phi(\lambda;x_1,\ldots,x_{m-1})$ is a solution to Equation (\ref{nod}).

If $s = 0$ then the greatest common divisor $d^{(m)}(\lambda,x_1,\ldots,x_m)$ does not generate the
solutions of System (\ref{newsys}) from the set $\Phi_{m-1}$, and, hence, solutions of System
(\ref{scsyst}).

If at least one of the functions $\widetilde {\widetilde{f}}^{(m)}_j(\lambda,x_1,\ldots,x_m)$ $(j =
1,\ldots,n)$ is distinct from zero at the zero point (this can occur only in the case when the
degree of the polynomial $d^{(m)}(\lambda,x_1,\ldots,x_m)$ coincides with the degree of one of the
polynomials $\widetilde{f}^{(m)}_j(\lambda,x_1,\ldots,x_m)$ $(j = 1,\ldots,n)$) then the process of
the construction of the set $\Phi$ is finished. Thus $\Phi$ coincides with $\Phi_{m-1}$, and
$\Phi_j$ ($j = 0,\ldots,m-2$) are the empty sets. Otherwise (i.e., when the degree $s$ of greatest
general divisor of the polynomials $\widetilde {f}^{(m)}_j(\lambda,x_1,\ldots,x_m)$ $(j =
1,\ldots,n)$ is strictly less the degrees of each of these polynomials) we pass to the
consideration of System (\ref{syx}).

Let us consider the system of the equations
 \begin{equation}\label{resultant}
 f_j^{(m-1)}(\lambda,x_1,\ldots,x_{m-1}) = 0 \quad (j=1,\ldots,n_{m-1}),
 \end{equation}
whose left hand parts are the full system of the resultants (see, for example, \cite{BM, KVZRS,
Leng}) for the polynomials standing in the left hand parts of System (\ref{syx}).

System (\ref{resultant}) is analogous to the initial system (\ref{scsyst}), however, its left hand
sides depend on the smaller number of the variables (namely from $\lambda,x_1,\ldots,x_{m-1}$).
Thus if $n = 2$, $n_{m-1} = 1$ and if $n > 2$, the number $n_{m-1}$ is greater the number $n$. The
following simple statement (see, for example, \cite {KVZRS}) will be used below.

{\bf Lemma 1.} {\it System of equations} (\ref{syx}) {\it has small solutions if and only if System
of equations} (\ref{resultant}) {\it has small solutions. More precisely, if $x(\lambda) =
(\xi_1(\lambda),\ldots,\xi_{m-1}(\lambda)$, $\xi_m(\lambda))$ is a small solution of System}
(\ref{syx}), {\it then $\widetilde{x}(\lambda) = (\xi_1(\lambda),\ldots,\xi_{m-1}(\lambda))$ is a
small solution of System} (\ref{resultant}). {\it Vice versa, if $\widetilde{x}(\lambda) =
(\xi_1(\lambda),\ldots,\xi_{m-1}(\lambda))$ is a small solution of System} (\ref{resultant}), {\it
then there is a finite (not equal to zero) number of continuous at zero and turning into zero at
zero functions $\xi_m(\lambda)$ for which $x(\lambda) =
(\xi_1(\lambda),\ldots,\xi_{m-1}(\lambda),\xi_m (\lambda))$ is a small solution of System}
(\ref{syx}).

One can apply to System (\ref{resultant}) the same argument as was applied to System
(\ref{scsyst}). Namely, reducing the left hand parts of System (\ref{resultant}) by the greatest
degrees of the parameter $\lambda$, then implementing the suitable linear substitution of the
unknowns and applying the Weierstrass preparation theorem, we see that System (\ref{resultant}) is
equivalent to the system of the algebraic equations with respect to the unknown $x_{m-1}$:
 \begin{equation}\label{newnewsys}
 \widetilde {f}^{(m-1)}_j(\lambda,x_1,\ldots,x_{m-1}) = 0 \qquad (j = 1,\ldots,n_{m-1}).
 \end{equation}

Let us denote by $d^{(m-1)}(\lambda,x_1,\ldots,x_{m-1})$ the greatest common divisor of the
polynomials $\widetilde{f}^{(m-1)}_j(\lambda $, $x_1,\ldots,x_{m-1}) \ (j = 1,\ldots,n_{m-1})$.
Then
 \begin{equation}\label{newnodd}
 \begin{array}{c}
 \widetilde{f}^{(m-1)}_j(\lambda,x_1,\ldots,x_{m-1}) = d^{(m-1)}(\lambda,x_1,\ldots,x_{m-1}) \cdot
 \widetilde{\widetilde{f}}^{(m-1)}_j(\lambda,x_1,\ldots,x_{m-1}) \\ (j = 1,\ldots,n_{m-1}).
 \end{array}
 \end{equation}
and system of equations (\ref{resultant}) is equivalent to the set of one algebraic equation of the
unknown $x_{m-1}$
 \begin{equation}\label{newnod}
 d^{(m-1)}(\lambda,x_1,\ldots,x_{m-1}) = 0,
 \end{equation}
and the system of algebraic equations of the unknown $x_{m-1}$
 \begin{equation}\label{newsix}
 \widetilde{\widetilde{f}}^{(m-1)}_j(\lambda,x_1,\ldots,x_{m-1}) = 0 \qquad (j = 1,\ldots,n_{m-1}).
 \end{equation}

If the degree $s_{m-1}$ of the greatest common divisor $d^{(m-1)}(\lambda,x_1,\ldots,x_{m-1})$ is
positive then equation (\ref{newnod}) at any enough small $\lambda,x_1,\ldots,x_{m-2}$ has
$s_{m-1}$ small solutions $x_{m-1} = \xi(\lambda,x_1,\ldots,x_{m-2})$. Thus these solutions will
define the elements of the set $\Phi_{m-2}$, that is the solutions of System (\ref {resultant})
 $$x(\lambda) = \phi(\lambda;x_1(\lambda),\ldots,x_{m-2}(\lambda)),$$
where $x_i(\lambda)$ ($i = 1,\ldots,m - 2$) are arbitrary free parameters. If $s_{m-1} = 0$ then
the greatest common divisor $d^{(m-1)}(\lambda,x_1,\ldots,x_{m-1})$ does not generate solutions of
System (\ref{newnewsys}), hence, it does not generate the solutions of System (\ref{resultant}).

If at least one of the functions $\widetilde{\widetilde{f}}^{(m-1)}_j(\lambda,x_1,\ldots,x_{m-1})$
$(j = 1,\ldots,n_{m-1})$ is distinct from zero at zero then system of equations (\ref{newsix}) has
no small solutions. One can construct solutions of system (\ref{scsyst}) by means of solutions of
equations (\ref{newnod}) and Lemma 1. The set of these solutions forms the set $\Phi_{m-2}$. Thus
$\Phi = \Phi_{m-2} \cup \Phi_{m-1}$ and the process of the construction of the set $\Phi$ is
finished. Otherwise we pass to the consideration of System (\ref{newsix}).

If we apply to System (\ref{newsix}) the argument which was applied to System (\ref{syx}) we
construct the set $\Phi_{m-3}$. Other sets $\Phi_k$ ($k = 0,\ldots,m - 4$) are constructed
similarly. After reduction by suitable degrees of parameter $\lambda$, implementation of the linear
change of variables and applying the Weierstrass preparation theorem to every system
 \begin{equation}\label {resultant-k}
 f_j^{(k)}(\lambda,x_1,\ldots,x_k) = 0 \quad (j = 1,\ldots,n_k)
 \end{equation}
we obtain the system
 \begin{equation}\label{newnewsys-k}
 \widetilde{f}^{(k)}_j(\lambda,x_1,\ldots,x_k) = 0 \qquad (j = 1,\ldots,n_k);
 \end{equation}
of the algebraic equations of the unknown $x_k$. Then the greatest common divisor
$d^{(k)}(\lambda,x_1,\ldots,x_{k})$ of the polynomials, that stand in the left hand parts of
equations of this system is defined. Finally, we construct the system
 \begin{equation}\label{newsix-k}
 \widetilde{\widetilde{f}}^{(k)}_j(\lambda,x_1,\ldots,x_k) = 0 \qquad (j = 1,\ldots,n_k).
 \end{equation}
At this moment the set $\Phi_{k-1}$ is nonempty if and only if the degree of the polynomial
$d^{(k)}(\lambda,x_1,\ldots,x_k)$ is positive; the set $\Phi_0 \cup \ldots \cup \Phi_{k-2}$ is
nonempty if and only if the degrees of all polynomials which stand in the left hand part of System
(\ref{newsix-k}) are positive or the left hand parts of the equations of System (\ref{resultant-k})
turn into zero at the zero values of the arguments.

The process of construction of the set $\Phi$ described above leads to a chain of sets of the
equations and systems of the equations. If the process does not break at some intermediate step,
then this chain can be present as:
 \begin{equation}\label{Lscheme}
 \begin{array}{l}
 \qquad \qquad \qquad \qquad f^{(m)}_j(\lambda,x_1,\ldots,x_m) = 0 \ \ (j = 1,\ldots,n) \\
 \qquad \qquad \qquad \qquad \swarrow \qquad \searrow \\
 d^{(m)}(\lambda,x_1,\ldots,x_m) = 0 \quad f_j^{(m-1)}(\lambda,x_1,\ldots,x_{m-1}) = 0 \ \ (j = 1,\ldots,n_{m-1}) \\
 \quad \qquad \qquad \qquad \qquad \qquad \swarrow \qquad \searrow \\
 \quad d^{(m-1)}(\lambda,x_1,\ldots,x_{m-1}) = 0 \quad f_j^{(m-2)}(\lambda,x_1,\ldots,x_{m-2}) = 0 \ \
 (j = 1,\ldots,n_{m-2}) \\
 \qquad \qquad \qquad \qquad \qquad \qquad \qquad \swarrow \qquad \searrow \\
 \qquad \ldots \ldots \ldots \ldots \ldots \ldots \ldots \ldots \ldots \ldots \ldots \ldots \ldots \ldots
 \ldots \ldots \ldots \ldots \ldots \ldots \ldots \ldots \ldots \\
 \qquad \qquad \qquad \qquad \qquad \qquad \quad \swarrow \qquad \searrow \\
 \quad \quad d^{(k+1)}(\lambda,x_1,\ldots,x_{k+1}) = 0 \quad f_j^{(k)}(\lambda,x_1,\ldots,x_k) = 0 \
 (j=1,\ldots,n_k) \\
 \quad \qquad \qquad \qquad \qquad \qquad \qquad \swarrow \qquad \searrow \\
 \quad \quad \quad d^{(k)}(\lambda,x_1,\ldots,x_k) = 0 \quad f_j^{(k-1)}(\lambda,x_1,\ldots,x_{k-1}) = 0 \
 (j=1,\ldots,n_{k-1}) \\
 \qquad \qquad \qquad \qquad \qquad \qquad \swarrow \qquad \searrow \\
 \qquad \quad \quad \ldots \ldots \ldots \ldots \ldots \ldots \ldots \ldots \ldots \ldots \ldots \ldots
 \ldots \ldots \ldots \ldots \ldots \ldots \ldots \ldots \ldots \\
 \qquad \qquad \qquad \qquad \qquad \qquad \swarrow \qquad \searrow \\
 \quad \quad \quad \quad d^{(3)}(\lambda,x_1,x_2,x_3) = 0 \quad f_j^{(2)}(\lambda,x_1,x_2) = 0 \ (j=1,\ldots,n_2) \\
 \qquad \qquad \qquad \qquad \qquad \qquad \quad \swarrow \qquad \searrow \\
 \quad \quad \quad \quad \quad d^{(2)}(\lambda,x_1,x_2) = 0 \quad f_j^{(1)}(\lambda,x_1) = 0 \ (j=1,\ldots,n_1) \\
 \qquad \qquad \qquad \qquad \qquad \qquad \swarrow \\
 \quad \quad \quad \quad \quad \quad d^{(1)}(\lambda,x_1) = 0.
 \end{array}
 \end{equation}

Let us write out from (\ref{Lscheme}) the equations participating in the construction of the set
$\Phi_k$ ($k = 0,\ldots,m - 1$). We arrange them in the form
 \begin{equation}\label{chain}
 \begin{array}{l}
 f^{(m)}_j(\lambda,x_1,\ldots,x_m) = 0 \quad (j=1,\ldots,n) \\ \qquad \qquad \qquad \qquad \qquad
 \quad \nwarrow \\ \qquad f_j^{(m-1)}(\lambda,x_1,\ldots,x_{m-1}) = 0 \quad (j=1,\ldots,n_{m-1}) \\
 \qquad \qquad \qquad \quad \qquad \qquad \quad \nwarrow \\ \qquad \qquad \qquad \qquad \qquad
 \ldots \ldots \ldots \ldots\ldots \\ \qquad \qquad \qquad \qquad \qquad \qquad \quad \nwarrow \\
 \qquad \qquad \qquad f_j^{(k+1)}(\lambda,x_1,\ldots,x_{k+1}) = 0 \quad (j=1,\ldots,n_{k+1}) \\
 \qquad \qquad \qquad \quad \qquad \qquad \nearrow \\ \qquad d^{(k+1)}(\lambda,x_1,\ldots,x_{k+1}) = 0
 \quad (k=0,\ldots,m-1). \end{array}
 \end{equation}

Diagram (\ref{chain}) shows that the solutions of the equation
 $$d^{(k+1)}(\lambda,x_1,\ldots,x_{k+1}) = 0$$
define $(k+1)$-th components
 $$x_{k+1}(\lambda) = \psi_{k+1}(x_1(\lambda),\ldots,x_k (\lambda))$$
of the required solutions $x_1(\lambda),\ldots,x_k(\lambda)$. Further we pass to the system of the
equations
 \begin{equation}\label{v1}
 f_j^{(k+1)}(\lambda,x_1,\ldots,x_{k+1}) = 0 \quad (j=1,\ldots,n_{k+1}).
\end{equation}
The solutions of System (\ref{v1}) define $(k+2)$-th components
 $$x_{k+2}(\lambda) = x_{k+2}(x_1(\lambda),\ldots,x_k(\lambda),x_{k+1}(\lambda))$$
of the required solutions. Moving <<upwards>> we pass to the next system of equations and etc.
Finally we <<reach>> the last system of equations
 \begin{equation}\label{vv}
 f^{(m)}_j (\lambda, x_1, \ldots, x_m) = 0 \quad (j=1, \ldots, n).
 \end{equation}
The solutions of System (\ref{vv}) define the last components
 $$x_m (\lambda) = x_m(x_1(\lambda),\ldots,x_k(\lambda),x_{k+1}(\lambda),\ldots,x_{m-1}(\lambda))$$
of the required solutions.

As a consequence of (\ref{Lscheme}) and (\ref{chain}) one obtains the following \cite{KVZRS,VT}:

{\bf Theorem 1.} {\it In the complex case System} (\ref{scsyst}) {\it has a finite number of small
solutions if and only if the degrees of the polynomials $d^{(i)}(\lambda,x_1,\ldots,x_i)$ ($i =
2,\ldots,m$) are equal to zero. Thus, if the degree of the polynomial $d^{(1)}(\lambda,x_1)$ is
equal to zero then System} (\ref{scsyst}) {\it has no small solutions. If this degree is positive,
then System} (\ref{scsyst}) {\it has a finite number of small solutions.}

The simple examples show that in the real case the analogue of Theorem 1 is false.

{\bf Analysis of the R\"uckert--Lefschetz scheme.} The described scheme of research of the implicit
functions defined by System (\ref{scsyst}) allows to describe the general structure of the small
solutions $x = x(\lambda)$ of this system. One naturally comes to the following question: Is it
possible, using the R\"uckert--Lefschetz scheme, to construct implicit functions (i.e. elements of
the set $\Phi$) determined by System (\ref{scsyst})? The matter is that calculation for analytic
functions is usually realized through their expansion into the series of their variables (such
calculation is usually called {\it approximate}). Moreover, actually only the first coefficients
are used in the calculation. Therefore there appears a new question: Are the first coefficients in
Taylor expansions of solutions to System (\ref{scsyst}) determined by the first coefficients in the
Taylor expansions of the left hand parts of this system in Taylor series?

Recall that calculation in the R\"uckert--Lefschetz scheme described above is really realized in
the rings of analytic at zero functions of the variables $(\lambda,x_1,\ldots,x_{m-1},x_m)$,
$(\lambda,x_1,\ldots,x_{m-1})$, \ldots , $(\lambda,x_1)$, $\lambda$. However, dealing with concrete
systems of equations we must operate only with a finite number of coefficients in the corresponding
Taylor expansions of solutions and left hand parts of the system under consideration, or --- as
accepted to speak in the problem under consideration --- to within the members of the higher order.
At the first sight it seems that one can probably carry out similar <<approximate>> calculation for
the R\"uckert--Lefschetz scheme described above. However, it is not so. More precisely, in the
process of such calculation for concrete systems one can meet situations when the first
coefficients of the Taylor expansions of solutions are not determined by the first coefficients of
the Taylor expansions of the left hand parts of the system under consideration.

More exactly: {\it if calculation due to the R\"uckert--Lefschetz scheme is realized in the
framework of calculation exploiting only a finite number of coefficients (not in the rings of
analytic functions or the corresponding formal power series!) then it, generally speaking, does not
allow to define a number of these implicit functions and jets of these solutions}. Thus, the
R\"uckert--Lefschetz scheme has different properties in the framework of calculation in the rings
of analytic functions and in the framework of calculation with the jets of the left hand sides of
equations of System (\ref{scsyst}) and jets of its small solutions.

Of course, in the simplest case $m = n = 1$ we do not meet this problem. The Newton diagram method
states that the answer to these questions is positive for simple solutions (a solution $x(\lambda)$
of the scalar equation $f(\lambda,x) = 0$ is simple if $f'_x(\lambda,x(\lambda))$ is not zero) and
negative for not simple solutions.

Although we are interested in the cases when $m = n$, under the realization of the
R\"uckert--Lefschetz scheme we arrive at the systems with $m\neq n$ (more precisely, with $m < n$).
Let us consider one of such cases: $m = 1$, $n = 2$. System (\ref{scsyst}) in this case has the
form
 \begin{equation}\label{scsyst222}
 \left\{
 \begin{array}{l}
 f_1(\lambda,x_1) = 0, \\
 f_2(\lambda,x_1) = 0.
 \end{array}
 \right.
 \end{equation}
The scheme described above in this case leads to calculation of the resultant $f_{12}$ of the left
hand parts of System (\ref{scsyst222}) and to analysis of the equation
 \begin{equation}\label{scsyst111}
 f_{12}(\lambda) = 0.
 \end{equation}
The equations of System (\ref{scsyst222}) have a common solution if and only if the resultant
$f_{12}(\lambda) = 0$. However, this equality is determined by an infinite number of the
corresponding coefficients of $f_{12}(\lambda)$, and calculation of the latter ones requires the
knowledge of an infinite number of coefficients of the left hand parts of System (\ref{scsyst222}).
Thus, even for this simplest overdetermined system its solvability depends not only on the first
coefficients in Taylor expansions of the left hand parts of the system underconsideration.
Evidently, an analogous statement holds for an arbitrary overdetermined system.

The next simple case is $m=2$, $n=1$. System (\ref{scsyst}) in this case has the form
 \begin{equation}\label{scsyst333}
 f_1(\lambda,x_1,x_2) = 0.
 \end{equation}
and, without the loss of generality, one can assume that the left hand side $f_1(\lambda,x_1,x_2)$
of this equation is a polynomial with respect to $x_2$. In spite of the simplicity of this equation
the analysis of its solutions is seriously difficult and requires the usage of Singularity Theory.
However, application of the R\"uckert--Lefschetz scheme, in the main (for us) case $m = n$, leads
only to determined systems (if $m = n = 2$) or overdetermined systems (if $m = n > 2$). Thus, in
the analysis of System (\ref{scsyst333}) we need more complicated underdetermined systems.

Let us consider the case: $m = n = 2$. In this case System (\ref{scsyst}) has the form
 \begin{equation}\label{scsyst2}
 \left\{
 \begin{array}{l}
 f_1(\lambda,x_1,x_2) = 0, \\
 f_2(\lambda,x_1,x_2) = 0.
 \end{array}
 \right.
 \end{equation}
The scheme described above in this case leads to calculation of the resultant $f_{12}$ of the left
hand parts of System (\ref{scsyst2}) and to analysis of the equation
 \begin{equation}\label{scsyst1}
 f_{12}(\lambda,x_1) = 0.
 \end{equation}

From description of the R\"uckert--Lefschetz scheme it follows that the first coefficients in
Taylor expansions of the left hand sides of System (\ref{scsyst2}) determine the first coefficients
in Taylor expansion of the left hand sides of System (\ref{scsyst1}). Then, in generic cases, in
order to analyze System (\ref{scsyst1}) it is possible to apply the Newton diagram method. This
allows one to define, generally speaking, the first members of all solutions of System
(\ref{scsyst1}).

Further, substituting these approximate solutions $x_1(\lambda)$ into the equations of System
(\ref{scsyst2}) we receive a system of compatible equations for the definition of the second
components $x_2(\lambda)$ of the solutions of System (\ref{scsyst2}):
 \begin{equation}\label{scsyst3}
 \left\{
 \begin{array}{l}
 f_1(\lambda,x_1(\lambda),x_2) = 0, \\
 f_2(\lambda,x_1(\lambda),x_2) = 0.
 \end{array}
 \right.
 \end{equation}
This system is similar to System (\ref{scsyst222}) (with the unknown $x_2$ instead of $x_1$),
however, now we know that this system is solvable. Applying the Newton diagram method to each
equation of System (\ref{scsyst3}) one can construct the jets of all solutions to each equation of
System (\ref{scsyst3}). If there exists a unique common jet of solutions to equations of System
(\ref{scsyst3}) then this jet is a jet of a common solution to both equations of System
(\ref{scsyst3}). In all other cases we can only state that System (\ref{scsyst3}) is solvable but
does not determine jets of common solutions to System (\ref{scsyst3}).

Really, if ${\cal S}_i = \{x_2^{i,\sigma}(\lambda): \ \sigma = 1,\ldots,s_i\}$, $i = 1,2$ are the
sets of solutions of System (\ref{scsyst3}), then the set of solutions of System (\ref{scsyst3})
coincides with the set ${\cal S}_1 \cap {\cal S}_2$. However, the previous argument shows that we
can deal only with jets of the corresponding solutions. These jets form the new sets
$\widetilde{{\cal S}}_i = \{\widetilde{x}_2^{i,\sigma} (\lambda): \ \sigma = 1,\ldots,s_i\}$. In
the case under consideration the intersection $\widetilde{{\cal S}}_1 \cap \widetilde{{\cal S}}_2$
contains at least two elements. And in this case it is impossible to determine which of them really
determines common solutions to System (\ref{scsyst3}) and which no.

Thus, using the R\"uckert--Lefschetz scheme it is possible to determine the first coefficients in
Taylor expansions of solutions to System (\ref{scsyst2}) if for each solution $x_1(\lambda)$ to
System (\ref{scsyst1}) there exists a unique common jet of solutions to equations of System
(\ref{scsyst3}). This common jets determine the second components $x_2(\lambda)$ (for each
$x_1(\lambda)$) of solution to System (\ref{scsyst2}).

Here, it must be emphasized that all the examples of concrete systems with two equations and two
unknowns presented in monograph \cite{VT} are covered by the unique, pointed out above, case when
the R\"uckert--Lefschetz scheme allows us to construct jets of solutions.

Now we pass to the case, when $m = n > 2$. The corresponding system has the form
 \begin{equation}\label{scsystn}
 \left\{
 \begin{array}{l}
 f_1(\lambda,x_1,\ldots,x_{n-1},x_n) = 0, \\ \hdotsfor 1 \\
 f_n(\lambda,x_1,\ldots,x_{n-1},x_n) = 0,
 \end{array}
 \right.
 \end{equation}
and, without the loss of generality, we can assume that each $f_j(\lambda,x_1,\ldots,x_{n-1},x_n)$,
$j = 1,\ldots,n$, is a polynomial of positive degree with respect to $x_n$. In the framework of
calculation exploiting the first coefficients in expansions of these polynomials we ought to
consider only the case when all members of resultant system to these polynomials are nonzero. In
addition, the number of members in the obtained resultant system is more than $n$. Therefore, the
corresponding system of equations
 \begin{equation}\label{scsystnr}
 \left\{
 \begin{array}{l}
 \widetilde{f}_1(\lambda,x_1,\ldots,x_{n-1}) = 0, \\ \hdotsfor 1 \\
 \widetilde{f}_{n_1}(\lambda,x_1,\ldots,x_{n-1}) = 0,
 \end{array}
 \right.
 \end{equation}
($n_1 > n - 1$) is overdetermined. We can chose $n - 1$ equations among equations of this resultant
system and, in a generic case, find small solutions $(x_1(\lambda),\ldots,x_{n-1}(\lambda))$,
defined by this system of $n - 1$ equations with $n - 1$ unknowns. Further among these solutions we
must gather those satisfying other equations of System (\ref{scsystnr}). However, this can be done
only if we use infinite number of coefficients in Taylor expansions of the left hand sides of these
equations. The latter is not possible in the framework of calculation with the first coefficients.

Thus, the R\"uckert--Lefschetz scheme, in the framework of calculation with the first coefficients,
does not allow, generally speaking, to determine the first coefficients of the Taylor expansions of
solutions to System (\ref{scsyst}) (even for a rough systems). To give the exact description of
this fact we need a new definition.

Let ${\frak M}$ be a class of finite systems of type (\ref{scsyst}) with analytical left hand
parts. We say that some scheme (algorithm) ${\frak S}$ of investigation of solutions of systems
from ${\frak M}$ is {\it effective} if this scheme allows to define jets of all solutions of a
system from ${\frak M}$ using only a finite number of the first coefficients in the Taylor
expansions of the left hand parts of the system under consideration. It is evident, that the
R\"uckert--Lefschetz scheme is non effective in the class ${\frak M}$ of finite systems of type
(\ref{scsyst}) if in this class there exists a system without the property of roughness. So, the
R\"uckert--Lefschetz scheme can be effective only in the case when the class ${\frak M}$ contains
only rough systems. However, the above stated argument proves the following statement.

{\bf Theorem 2.} {\it The R\"uckert--Lefschetz scheme of construction of small solutions of system
of type (\ref{scsyst}) is not effective for $m = n > 1$ for the class of rough systems.}

\smallskip

Let us remind (see for example \cite {KVZRS}) that the Newton diagram method of investigation of
one scalar equation $f(\lambda,x) = 0$ with an analytical left hand part is effective.

\smallskip

{\bf Refinement of the R\"uckert--Lefschetz scheme.} Below we give some standard complements to the
R\"uckert--Lefschetz scheme, although these complements lie outside the theme of our main results.

Let us consider the case, when $k = 0$. We present the polynomial $q^{(1)}(\lambda,x_1) =
d^{(1)}(\lambda,x_1)$ in the form of the product of prime multipliers over the ring
$K[\lambda,x_1]$ of analytic at zero functions. Let $p^{(1)}(\lambda,x_1)$ is one of the prime
multipliers of the polynomial $q^{(1)}(\lambda,x_1)$. Then each solution $x_1(\lambda)$ of the
equation
 \begin{equation}\label{p1}
 p^{(1)}(\lambda,x_1) = 0
 \end{equation}
is the first component of an element in the set $\Phi_0$.

To define the second components $x_2(\lambda)$ of elements in the set $\Phi_0$ whose first
components are solutions of Equation (\ref{p1}) we consider the following system from
(\ref{Lscheme})
 \begin{equation}\label{z1}
 f_j^{(2)}(\lambda,x_1,x_2) = 0 \quad (j = 1,\ldots,n_2).
 \end{equation}
The left hand parts of the equations in this system are polynomials in $x_2$ with coefficients from
the ring $K[\lambda,x_1]$ of analytic at zero functions turning into zero at zero.

According to the Weierstrass preparation theorem let us replace coefficients of the polynomials
standing in the left hand part of System (\ref{z1}) with their remainders from division of them by
the prime polynomial $p^{(1)}(\lambda,x_1)$. As the result of such replacement System (\ref{z1})
transforms into the system
 \begin{equation}\label{n1}
 \widehat{f}_j^{(2)}(\lambda,x_1,x_2) = 0 \quad (j = 1,\ldots,n_2),
 \end{equation}
where $\widehat{f}_j^{(2)}(\lambda,x_1,x_2)$ $(j = 1,\ldots,n_2)$ are polynomials in $x_2$ whose
coefficients are polynomials in $x_1$ and whose degrees are less than the degree of the polynomial
$p^{(1)}(\lambda,x_1)$.

Since $p^{(1)}(\lambda, x_1)$ is a prime polynomial it is possible to consider System (\ref{n1}) as
a system of algebraic equations with the unknown $x_2$ in the field $K(\lambda,x_1)$ which was
obtained from the field $K(\lambda)$ by adding an algebraic element $x_1$, where $x_1$ is a
solution of Equation (\ref{p1}). Since the concept of the greatest common divisor is defined in the
ring of polynomials over the field, System (\ref{n1}) is equivalent to the single equation
 \begin{equation}\label{o1}
 q^{(2)}(\lambda,x_1,x_2) = 0,
 \end{equation}
where $q^{(2)}(\lambda,x_1,x_2)$ is the greatest common divisor of the polynomials standing in the
left hand part of System (\ref{n1}). Thus, to define the second component $x_2(\lambda)$ of
elements in the set $\Phi_0$ whose first components $x_1 (\lambda)$ are solutions of Equation
(\ref{p1}), it is enough to find solutions of the algebraic equation (\ref{o1}).

% with coefficients from the field $K(\lambda,\xi_1)$, where $\xi_1$ is an algebraic over $K(\lambda)$
%element coinciding with a root of Equation (\ref{p1}).

Let $p^{(2)}(\lambda,x_1,x_2)$ be one of the prime multipliers of the polynomial
$q^{(2)}(\lambda,x_1,x_2)$. We find the third components $x_3(\lambda)$ of elements in the set
$\Phi_0$ whose first components $x_1(\lambda)$ are solutions of Equations (\ref{p1}) and second
components are solutions of the equation
 \begin{equation}\label{p2}
 p^{(2)}(\lambda,x_1,x_2) = 0.
 \end{equation}
To this end we consider the following system from (\ref{Lscheme})
 \begin{equation}\label{z2}
 f_j^{(3)}(\lambda,x_1,x_2,x_3) = 0 \quad (j = 1,\ldots,n_3).
 \end{equation}
The left hand parts of the equations of System (\ref{z2}) are polynomials in $x_3$ with
coefficients from the ring $K[\lambda,x_1,x_2]$ of analytic at zero functions turning into zero at
zero. Considering System (\ref{z2}) together with Equation (\ref{p2}) one can simplify System
(\ref{z2}). Firstly, in accordance with the Weierstrass preparation theorem, it is possible to
replace each coefficient of the left hand parts of System (\ref{z2}) with the remainder of its
division by the prime polynomial $p^{(2)}(\lambda,x_1,x_2)$. And secondly, it is possible to divide
the coefficients of the polynomials obtained by the prime polynomial $p^{(1)}(\lambda,x_1)$
according to the Weierstrass preparation theorem and replace them with the remainders from these
divisions. As the result System (\ref{z2}) transforms into the system
 \begin{equation}\label{n2}
 \widehat{f}_j^{(3)}(\lambda,x_1,x_2,x_3) = 0 \quad (j = 1,\ldots,n_3),
 \end{equation}
where $\widehat{f}_j^{(3)}(\lambda,x_1,x_2,x_3) = 0$ $(j = 1,\ldots,n_3)$ are polynomials in $x_3$
whose coefficients are polynomials in $x_2$ and whose degrees are less than the degree of the
polynomial $p^{(2)}(\lambda,x_1,x_2)$ and coefficients of these polynomials are polynomials in
$x_1$ whose degrees are less than the degree of the polynomial $p^{(1)}(\lambda,x_1)$.

Since $p^{(2)}(\lambda,x_1,x_2)$ is a prime polynomial it follows that System (\ref{n2}) can be
considered as a system of algebraic equations with the unknown $x_3$ in the field
$K(\lambda,x_1,x_2)$ which was obtained from the field $K(\lambda,x_1)$ by adding an algebraic
element $x_2$, where $x_2$ is a solution of Equation (\ref{p2}). In this case System (\ref{n2}) is
equivalent to the single equation
 \begin{equation}\label{o2}
 q^{(3)}(\lambda,x_1,x_2,x_3) = 0,
 \end{equation}
where $q^{(3)}(\lambda,x_1,x_2,x_3)$ is the greatest common divisor of the polynomials standing in
the left hand part of System (\ref{n2}). To determine the third component $x_3(\lambda)$ of
solutions $x(\lambda)$ of System (\ref{scsyst}) which first components are solutions of Equation
(\ref{p1}) and the second components are solutions of Equation (\ref{p2}) it is enough to find
solutions of the algebraic equation (\ref{o2}). Thus, to define the third component $x_3(\lambda)$
of elements in the set $\Phi_0$ whose first component $x_1 (\lambda)$ is a solution of Equation
(\ref{p1}) and whose second component $x_2(\lambda)$ is a solution of Equation (\ref{p2}), it is
enough to find solutions of the algebraic equation (\ref{o2}).

Continuing similarly we show that each component $x_j(\lambda)$ ($j = 1,\ldots,m$) of elements in
the set $\Phi_0$ will be defined by the algebraic equation
 \begin{equation}\label{pj}
 p^{(j)}(\lambda,x_1,\ldots,x_{j-1},x_j) = 0,
 \end{equation}
whose left hand part is a prime polynomial in variable $x_j$ with coefficients from the ring
$K[\lambda,x_1,\ldots,$ $x_{j-1}]$ of analytic at zero functions turning into zero at zero.
Moreover these coefficients are polynomials in variable $x_{j-1}$ with coefficients from the ring
$K[\lambda,x_1,\ldots,x_{j-2}]$ of analytic at zero functions turning into zero at zero. In turn,
the coefficients of these polynomials are polynomials in variable $x_{j-2}$ with coefficients from
the ring $K[\lambda,x_1,\ldots,x_{j-3}]$ of analytic at zero functions turning into zero at zero,
and etc.

Collecting equations (\ref{p1}), (\ref{p2}), (\ref{pj}) ($j = 3,\ldots,m$) we get that each element
in the set $\Phi_0$ is defined by a system of algebraic equations
 \begin{equation}\label{tt}
 \left\{\begin{array}{r}
 p^{(m)}(\lambda,x_1,\ldots,x_k,\ldots,x_m) = 0, \\
 \ldots\ldots\ldots\ldots\ldots\ldots\ldots\ldots\ldots \\
 p^{(k)}(\lambda,x_1,\ldots,x_k) = 0, \\
 \ldots\ldots\ldots\ldots\ldots\ldots\ldots \\
 p^{(2)}(\lambda,x_1,x_2) = 0, \\
 p^{(1)}(\lambda,x_1) = 0. \\
 \end{array}\right.
 \end{equation}
with coefficients, whose structure is described above.

%where each function $p^{(j)}(\lambda,x_1,\ldots,x_j)$ ($j = 1,\ldots,m$) is a prime polynomial of
%variable $x_j$ with coefficients from the ring $K[\lambda,x_1,\ldots,x_{j-1}]$. In others words
%the first component $x_1(\lambda)$ of an element from set $\Phi_0$ can be considered as an element
%of algebraic expansions of the field of quotients $K(\lambda)$ from the ring $K[\lambda]$ of
%analytical functions at the zero, the second component $x_2(\lambda)$ of this element can be
%considered as an element of the algebraic expansion of the field of the quotients $K(\lambda,x_1)$
%of the ring $K[\lambda,x_1]$ of analytical functions at zero, ..., at last, the last component
%$x_m(\lambda)$ of this element can be considered as an element of the algebraic expansion of the
%field of the quotients $K(\lambda,x_1,\ldots,x_{m-1})$ from the ring $K[\lambda,x_1,\ldots,x_{m-1}]$
%of analytical functions at the zero.

It is known (see, for example, \cite {L,VW}) that a finite number of the consecutive algebraic
expansions of the field of the quotients $K(\lambda)$ is equivalent to a simple algebraic
expansion. Namely, there are complex numbers $a_i$ ($i = 1,\ldots,m$) for which: (i) the function
$\eta(\lambda) = \sum\limits_{i=1}^m a_i \cdot x_i (\lambda)$ satisfies an algebraic equation
$\psi(\lambda,\eta) = 0$ with coefficients from the field $K(\lambda)$ and (ii) each function $x_j
(\lambda)$ ($j = 1,\ldots,m$) lies in the field $K(\lambda,\eta)$, i.e. has the form $x_j(\lambda)
= c_{j1}(\lambda)+c_{j2}(\lambda)\eta(\lambda)+\ldots+c_{js}(\lambda)\eta^{s-1}(\lambda)$, where
$c_{j\sigma}(\lambda)$ ($j = 1,\ldots,m$, $\sigma = 1,\ldots,s$) are functions from $K(\lambda)$,
and $s$ is the degree of the equation $\psi(\lambda, \eta) = 0$. Thus, System (\ref{tt}) is
equivalent to a single algebraic equation and the components of the corresponding element in the
set $\Phi_0$ are defined by the roots of this algebraic equation:
 \begin{equation}\label{ttt}
 \left\{\begin{array}{l}
 \psi(\lambda,\eta) = 0, \\
 x_j (\lambda) = c_{j1}(\lambda)+c_{j2}(\lambda)\eta(\lambda)+ \ldots + c_{js}(\lambda)\eta^{s-1}(\lambda)
 \quad (j = 1,\ldots,m).
 \end{array}\right.
 \end{equation}

Let us remind that the field of fractions $K(\lambda)$ can be presented as
 \begin{equation}\label{field}
 K(\lambda) = \left\{\lambda^\theta\sum_{i=0}^p v_i\lambda^i: \ v_i \in {\Bbb C}, \ v_0 \neq 0, \
 p \in {\Bbb N}, \ \theta \in {\Bbb Z} \right\}.
 \end{equation}
The field of fractions $K(\lambda)$ is not algebraically closed, however its algebraic closure
$K^*(\lambda)$ can be easily described
 \begin{equation}\label{closefield}
 K^*(\lambda) = \bigcup_{r=1}^\infty K \big (\lambda^\frac1r\big),
 \end{equation}
where $K \big(\lambda^\frac1r \big)$ is the field of fractions of the rings $K \big[\lambda^\frac1r
\big]$ of analytic functions.

From the aforesaid it follows that each element in the set $\Phi_0$ is an analytic function in the
parameter $\lambda$ or an analytic function in $\lambda^\frac1r$.

Let us notice also that each system (\ref{tt}) defines one or several elements in the set $\Phi_0$.

Let us pass to the case when $k>0$. The basic argument which was exploited in the analysis of the
set $\Phi_0$ is saved for the analysis of the sets $\Phi_k$.

Let us present the polynomial $q^{(k+1)}(\lambda,x_1,\ldots,x_k,x_{k+1}) =
d^{(k+1)}(\lambda,x_1,\ldots,x_k,x_{k+1})$ in the form of the product of prime polynomials over the
ring $K[\lambda,x_1,\ldots,x_k]$ of analytic functions at zero.

Let $p^{(k+1)}(\lambda,x_1,\ldots,x_k,x_{k+1})$ be one of the prime multipliers of the polynomial
$q^{(k+1)}(\lambda, x_1,\ldots,$ $x_k,x_{k+1})$.

Let us choose arbitrary $x_1(\lambda),\ldots,x_k(\lambda)$. Then each solution $x_{k+1}$ of the
equation
 \begin{equation}\label{pp1}
 p^{(k+1)}(\lambda,x_1,\ldots,x_k,x_{k+1}) = 0
 \end{equation}
is $(k+1)$-th component of an element in the set $\Phi_k$.

To determine $(k+2)$-th components of the elements in the set $\Phi_k$ whose first $k$ components
are arbitrary and $(k+1)$-th components $x_{k+1}(\lambda)$ are solutions of Equations (\ref{pp1})
we consider the system
\begin{equation}\label{zz1}
 f_j^{(k+2)}(\lambda,x_1,\ldots,x_k,x_{k+1},x_{k+2}) = 0 \quad (j = 1,\ldots,n_{k+2}).
 \end{equation}
Repeating the argument exploited above in the construction of the second components of solutions in
the set $\Phi_0$, we pass to the equivalent system
 \begin{equation}\label{nn1}
 \widehat{f}_j^{(k+2)}(\lambda,x_1,\ldots,x_k,x_{k+1},x_{k+2}) = 0 \quad (j = 1,\ldots,n_{k+2}),
 \end{equation}
where $\widehat{f}_j^{(k+2)}(\lambda,x_1,\ldots,x_k,x_{k+1},x_{k+2}) = 0$ $(j = 1,\ldots,n_{k+2})$
are polynomials in $x_{k+2}$ whose coefficients are polynomials in $x_{k+1}$ and degrees are less
than the degree of the polynomial $p^{(k+1)}(\lambda,x_1,\linebreak\ldots,$ $x_k,x_{k+1})$. Then we
pass to the equation
 \begin{equation}\label{oo1}
 q^{(k+2)}(\lambda,x_1,x_2,\ldots,x_k,x_{k+1},x_{k+2}) = 0,
 \end{equation}
where $q^{(2)}(\lambda,x_1,x_2)$ is the greatest common divisor of the polynomials standing in the
left hand part of System (\ref{nn1}). Thus, to determine $(k+2)$-th components $x_{k+2}(\lambda)$
of elements in the set $\Phi_k$ whose $(k+1)$-th components are solutions of Equation (\ref{pp1}),
it is enough to find the solutions of the algebraic equation (\ref{oo1}).

Continuing similarly, as well as in the case $k = 0$ we show that each component $x_j(\lambda)$ ($j
= k+1,\ldots,m$) of elements in the set $\Phi_k$ is defined by the algebraic equation
 \begin{equation}\label{ppj}
 p^{(j)}(\lambda,x_1,\ldots,x_k,x_{k+1},\ldots,x_{j-1},x_j) = 0,
 \end{equation}
whose left hand part is a prime polynomial in the variable $x_j$ with coefficients from the ring
$K[\lambda,x_1,\ldots,$ $x_{j-1}]$ of analytic functions at zero turning into zero at zero.
Moreover, in turn these coefficients are polynomials in the variable $x_{j-1}$ with coefficients
from the ring $K[\lambda,x_1,\ldots,x_{j-2}]$ of analytic functions at zero turning into zero at
zero, and in turn coefficients of these polynomials are polynomials in the variable $x_{j-2}$ with
coefficients from the ring $K[\lambda,x_1,\ldots,x_{j-3}]$ of analytic functions at zero turning
into zero at zero and etc.

Collecting equations (\ref{pp1}), (\ref{ppj}) ($j = k+2,\ldots,m$), we get that each element in the
set $\Phi_k$ is defined by the system
 \begin{equation}\label{ttttt}
 \left\{\begin{array}{r}
 p^{(m)}(\lambda,x_1,\ldots,x_k,x_{k+1},x_{k+2},\ldots,x_m) = 0, \\
 \ldots\ldots\ldots\ldots\ldots\ldots\ldots\ldots\ldots\ldots\ldots\ldots \\
 p^{(k+2)}(\lambda,x_1,\ldots,x_k,x_{k+1},x_{k+2}) = 0, \\
 p^{(k+1)} (\lambda,x_1,\ldots,x_k,x_{k+1}) = 0, \\
 \end{array}\right.
 \end{equation}
where $x_1,\ldots,x_k$ are free parameters and each function $p^{(j)}(\lambda,x_1,\ldots,x_j)$ ($j
= k+1,\ldots,m$) is a prime polynomial in the variable $x_j$ with coefficients from the ring
$K[\lambda,x_1,\ldots,x_{j-1}]$. In other words $(k+1)$-th component of an element in the set
$\Phi_k$ can be considered as an element from algebraic expansion of the field of fractions
$K(\lambda,x_1,\ldots,x_k)$ of the ring $K[\lambda,x_1,\ldots,x_k]$ of analytic functions at zero,
$(k+2)$-th component of this element can be considered as an element from algebraic expansion of
the field of fractions $K(\lambda,x_1,x_2,\ldots,x_k,x_{k+1})$ of the ring
$K[\lambda,x_1,x_1,\ldots,x_k,x_{k+1}]$ of analytic functions at zero, ..., at last, the last
component $x_m(\lambda)$ of this element can be considered as an element from algebraic expansions
of the field of fractions $K(\lambda,x_1,\ldots,x_{m-1})$ of the ring
$K[\lambda,x_1,\ldots,x_{m-1}]$ of analytic functions at zero.

As in the case $k = 0$ (see, for example, \cite{L,VW}), there are complex numbers $a_i$, ($i =
1,\ldots,m$) such that: (i) the function $\eta(\lambda) = \sum\limits_{i=k+1}^m a_i \cdot
x_i(\lambda)$ satisfies some algebraic equation $\psi(\lambda,x_1,\ldots,x_k,\eta) = 0$ with
coefficients from the field $K(\lambda,x_1,\ldots,x_k)$ and (ii) each function $x_j(\lambda)$ ($j =
k+1,\ldots,m$) lies in the field $K(\lambda,\eta)$, i. e. has the form $x_j(\lambda) =
c_{j1}(\lambda,x_1,\ldots,x_k) + c_{j2}(\lambda,x_1,\ldots,x_k)\eta(\lambda) + \ldots +
c_{js}(\lambda,x_1,\ldots,x_k)\eta^{s-1}(\lambda)$, where $c_{j\sigma}(\lambda,x_1,\ldots,x_k)$ ($j
= k+1,\ldots,m$, $\sigma = 1,\ldots,s$) are function from $K(\lambda,x_1,\ldots,x_k)$ and $s$ is
the degree of the equation $\psi(\lambda,x_1,\ldots,x_k$, $\eta) = 0$. Thus System (\ref{ttttt}) is
equivalent to a single algebraic equation and the components of elements in the set $\Phi_k$ are
defined by the root of this algebraic equation:
 \begin{equation}\label{tttt}
 \left\{\begin{array}{l}
 \psi(\lambda,x_1,\ldots,x_k,\eta) = 0, \\
 x_j(\lambda) = c_{j1}(\lambda,x_1,\ldots,x_k) + c_{j2}(\lambda,x_1,\ldots,x_k)\eta(\lambda) +
 \ldots + c_{js}(\lambda,x_1,\ldots,x_k)\eta^{s-1}(\lambda) \\ \phantom{00000000000000000000000000000000}
 \phantom{0000000000000000000000} (j = k + 1,\ldots,m).
 \end{array}\right.
 \end{equation}
Thus, each function $c_{j\sigma}(\lambda,x_1,\ldots,x_k)$ ($j = k+1,\ldots,m$, $\sigma =
1,\ldots,s$) can be presented in the form of fraction whose numerator and denominator are elements
of the ring $K[\lambda,x_1,\ldots,x_k]$ of analytic functions at zero. Therefore one has to choose
$x_1(\lambda),\ldots,x_k(\lambda)$ as the corresponding denominator of the functions
$c_{j\sigma}(\lambda,x_1,\ldots,x_k)$ not turning into zero at zero.

Let us notice that the analogues of formulas (\ref{field}) and (\ref{closefield}) do not exist in
the case $k>0$.

Again as well as in the case $k = 0$ each System (\ref{tttt}) defines one or several elements in
the set $\Phi_k$.

The argument presented above implies

{\bf Theorem 3.} {\it In the field of fractions $K (\lambda,x_1,\ldots,x_k) $ ($k=0,\ldots,m-1$) of
the ring $K[\lambda,x_1,\ldots,$ $x_k]$ of analytic functions at zero for each element in the set
$\Phi_k$ there is a prime equation $\psi(\lambda,x_1,\ldots,$ $x_k,\eta) = 0$ depending on free
parameters $\lambda$ and $x_1(\lambda),\ldots,x_k(\lambda)$, whose roots are defined by the
components $x_i(\lambda)$ $(i = k+1,\ldots,m-1)$ of this element. The components $x_i(\lambda)$ $(i
= k+1,\ldots,m-1)$ of an element in the set $\Phi_k$ depend on free parameters $\lambda$ and
$x_1(\lambda),\ldots,x_k(\lambda)$. Thus, the components of elements in the set $\Phi_0$ are
solutions of the systems whose equations have the form} (\ref {ttt}), {\it and the components of
elements in the sets $\Phi_k$ ($k = 1,\ldots,m-1$) are solutions of the systems of equations}
(\ref{tttt}).

It has been shown above that the R\"uckert--Lefschetz scheme is not an effective scheme for
construction of small solutions of System (\ref{scsyst}). One can see that the refinement of this
scheme is also "non effective" \ for calculation with coefficients of expansions of left hand sides
of the equations in System (\ref{scsyst}). However, below we show how to get the effective scheme
of the construction of some small solutions of System (\ref{scsyst}) by modifying the
R\"uckert--Lefschetz scheme.

{\bf Modified R\"uckert--Lefschetz scheme}. Applying the suitable change of variables and the
Weierstrass preparation theorem, we pass from System (\ref{scsyst}) to the consideration of the
equivalent system of algebraic (with respect to $x_n$) equations:
\begin{equation}\label{nnewsys}
\widetilde{f}^{(n)}_j(\lambda,x_1,\ldots,x_n) = 0 \qquad (j = 1,\ldots,n).
\end{equation}
Let us assume thus that $m = n$.

Here we shall exploit the set ${\frak D}_n$ of trees with $n$ vertexes. Let us remind that a tree
with $n$ vertexes is a coherent graph without simple cycles or, equivalently, a coherent graph with
$n$ vertexes and $n-1$ edges. The set ${\frak D}_n$ is finite; the number of its elements is equal
to $n^{n-2}$. A vertex of a tree is called multiple if it is an end vertex of more than one edge.
We denote by $\mu(D_n)$ the set of all multiple vertexes of a tree $D_n$.

Let $D_n$ be a tree from the set ${\frak D}_n$. Let us identify the vertexes of this tree with the
system of Equations (\ref{nnewsys}) (i.e. we enumerate vertexes of $D_n$ \ (1$ \le j \le n$) and
associate to $j$-th vertex of the tree $D_n$ the $j$-th equation of the system (\ref{nnewsys})).
Further, to each edge $\{j_1,j_2\}$ of the tree $D_n$ ($j_1$ and $j_2$ are numbers of the end
vertexes of the edge $\{j_1,j_2\}$) we associate the resultant of the left hand parts of $j_1$-th
and $j_2$-th equations of System (\ref{nnewsys}). As a result we get the system of $(n-1)$
equations with $n-1$ unknowns:
 \begin{equation}\label{v1a}
 \left\{\begin{array}{c}
 f^{(n-1)}_1 (\lambda,x_1,\ldots,x_{n-1}) = 0, \\
 \ldots \ldots \ldots \ldots \ldots \ldots \ldots \ldots \ldots \\
 f^{(n-1)}_{n-1}(\lambda,x_1,\ldots,x_{n-1}) = 0,
 \end{array}\right.
 \end{equation}
where $f^{(n-1)}_i(\lambda,x_1,\ldots,x_{n-1})$ ($i = 1,\ldots,n-1$) are the resultants
corresponding to the edges of $D$. We will assume that the left hand parts of this system are not
zero.

As it is known (see, for example, \cite{KVZRS}) $\widetilde{x}(\lambda) =
(\xi_1(\lambda),\ldots,\xi_{n-1}(\lambda))$ is a solution of System (\ref{v1a}), if $x(\lambda) =
(\xi_1(\lambda),\ldots,\xi_{n-1}(\lambda),\xi_n(\lambda))$ is a solution of System (\ref{scsyst}).
The opposite statement is not true. However, in some cases, it is possible to state that for a
given solution $\widetilde{x}(\lambda) = (\xi_1(\lambda),\ldots,\xi_{n-1}(\lambda))$ of System
(\ref{v1a}) there exists a unique solution $x(\lambda) =
(\xi_1(\lambda),\ldots,\xi_{n-1}(\lambda),\xi_n(\lambda))$ of System (\ref{scsyst}).

Let $\widetilde{x}(\lambda) = (\xi_1(\lambda),\ldots,\xi_{n-1}(\lambda))$ be a small solution of
System (\ref{v1a}). Let us consider the system
 \begin{equation}\label{v2}
 \left\{\begin{array}{c}f_1(\lambda,\xi_1(\lambda),\ldots,\xi_{n-1}(\lambda),x_n) = 0, \\
 \ldots \ldots \ldots \ldots \ldots \ldots \ldots \ldots \ldots \ldots \ldots \\
 f_n(\lambda,\xi_1(\lambda),\ldots,\xi_{n-1}(\lambda),x_n) = 0,\end{array} \right.
 \end{equation}
which is received from System (\ref{nnewsys}) by replacement of the unknowns $x_1,\ldots,x_{n-1}$
of this system by the components of the solution $\widetilde{x}(\lambda) =
(\xi_1(\lambda),\ldots,\xi_{n-1}(\lambda))$. We say that the solution $\widetilde{x}(\lambda) =
(\xi_1(\lambda),\ldots,\xi_{n-1}(\lambda))$ is {\it $D_n$-regular}, if System (\ref{v2}) has a
unique common simple solution $x_n = \xi_n(\lambda)$. According to this definition $D_n$-regular
solution $\widetilde{x}(\lambda) = (\xi_1(\lambda),\ldots,\xi_{n-1} (\lambda))$ uniquely defines
the solution $x(\lambda) = (\xi_1(\lambda),\ldots,\xi_{n-1}(\lambda),\xi_n(\lambda))$ of System
(\ref{nnewsys}).

At first sight it seems that the definition of $D_n$-regular solution of System (\ref{v1a}) is
senseless since this definition requires that the components of the solution
$\widetilde{x}(\lambda) = (\xi_1(\lambda),\ldots,\xi_{n-1}(\lambda))$ are the first components of
the corresponding solution of System (\ref{nnewsys}), i.e. this definition requires that the
solutions of System (\ref{nnewsys}) are defined by System (\ref{v1a}) that evidently is incorrect
in the general case. However, while in the general case it is impossible to construct all solutions
to System (\ref{nnewsys}) by means of solutions to System (\ref{v1a}), in some natural cases it is
possible to prove that a chosen solution to System (\ref{v1a}) is $D_n$-regular (certainly if it is
that) with the help of effective calculation (i.e. calculation using only a finite number of
coefficients in the expansion of the left hand parts of System (\ref{nnewsys})) and to construct
the missing component of the solution to System (\ref{nnewsys}).

The simple statement in this direction is

{\bf Lemma 2.} {\it If $\widetilde{x}(\lambda) = (\xi_1(\lambda),\ldots,\xi_{n-1}(\lambda))$ is a
small solution to systems} (\ref{v1a}) {\it and if each equation of System} (\ref{v2}) {\it with $j
\in \mu (D_n)$ has a unique solution then the solution $\widetilde{x}(\lambda) = \linebreak
(\xi_1(\lambda),\ldots,\xi_{n-1}(\lambda))$ to System} (\ref{v1a}) {\it is $D_n$-regular.}

Proof. Assume that $\widetilde{x}(\lambda) = (\xi_1(\lambda),\ldots,\xi_{n-1} (\lambda))$ is a
small solution to System (\ref{v1a}) and each equation of System (\ref{v2}) with $j \in \mu(D_n)$
has a unique small solution. In this case, if $j_1$ and $j_2$ are connected with an edge from $D_n$
then the corresponding equations
 \begin{equation}\label{ps}
 \begin{array}{l}
 f_{j_1}(\lambda,\xi_1(\lambda),\ldots,\xi_{n-1}(\lambda),x_n) = 0, \\
 f_{j_2}(\lambda,\xi_1(\lambda),\ldots,\xi_{n-1}(\lambda), x_n) = 0
 \end{array}
 \end{equation}
have a common (and unique) solution.

Let us consider equations of System (\ref{v2})
\begin{equation}\label{ps22}
 \begin{array}{l}
 f_{j_1}(\lambda,\xi_1(\lambda),\ldots,\xi_{n-1}(\lambda),x_n) = 0, \\
 f_{j_2}(\lambda,\xi_1(\lambda),\ldots,\xi_{n-1}(\lambda), x_n) = 0, \\
 \ldots \ldots \ldots \ldots \ldots \ldots \ldots \ldots \ldots \ldots \ldots\\
 f_{j_k}(\lambda,\xi_1(\lambda),\ldots,\xi_{n-1}(\lambda), x_n) = 0, \\
 \end{array}
 \end{equation}
where $j_1,j_2,\ldots,j_k \in \mu(D)$. Since $\mu(D_n)$ is a coherent subgraph of $D_n$, these
equations also have a common (and unique) solution $x_n = \xi_n(\lambda)$.

Now let us consider a pair of the equations of System (\ref{ps})
 \begin{equation}\label{ps2}
 \begin{array}{l}
 f_{j_1}(\lambda,\xi_1(\lambda),\ldots,\xi_{n-1}(\lambda),x_n) = 0, \\
 f_{j_2}(\lambda,\xi_1(\lambda),\ldots,\xi_{n-1}(\lambda),x_n) = 0,
 \end{array}
 \end{equation}
where $j_1 \in \mu(D_n)$, $j_2 \notin \mu(D_n)$, $j_1$ and $j_2$ are connected with an edge from
$D_n$. The first equation in this pair has a unique solution $\xi_n(\lambda)$. Therefore, both
equations have a common solution. So, $\xi_n(\lambda)$ is a solution of the second equation of this
system.

Thus, System (\ref{ps22}) has a unique common solution $\xi_n(\lambda)$ and, furthermore, System
(\ref{v2}) has a common (and unique) solution $\xi_n(\lambda)$. Hence System (\ref{nnewsys}) has a
common solution $\widetilde{x}(\lambda) = (\xi_1(\lambda),\ldots, \linebreak
\xi_{n-1}(\lambda),\xi_n(\lambda))$ and so this solution is $D_n$-regular.

Let us remind that if a solution $x_n(\lambda)$ of an equation of System (\ref{v1a}) is simple
(see, for example, \cite{KVZRS}), then beginning with a number $r$ all the subsequent coefficients
of the expansion of the simple solutions $x_n(\lambda)$ into the converging series in some
neighborhood of zero
 $$\begin{array}{c}
 x_n(\lambda) = \gamma_0 \lambda^\frac{\tau_0}{\tau} + \gamma_1
 \lambda^\frac{\tau_{1}}{\tau} + \ldots + \gamma_{l}
 \lambda^\frac{\tau_{l}}{\tau} + o(\lambda^\frac{\tau_{l}}{\tau}) \\[6pt]
 \quad (\gamma_{i} \in {\Bbb C}, \ \gamma_{0}\neq 0,\ \tau,
 \tau_{i} \in {\Bbb N}, \ i=0,\ldots,l)
 \end{array}$$
are defined by the
$$\alpha\gamma_{j} = \beta_{j} \quad (j\ge r),$$
where $\alpha$ is a constant. Hence the simple solution $x_n(\lambda)$ can be defined by a finite
number of coefficients of this expansion.

Lemma 2 is a special case of the following more general and obvious statement.

{\bf Lemma 3.} {\it Let $\widetilde {x}(\lambda) = (\xi_1(\lambda),\ldots,\xi_{n-1}(\lambda))$ be a
small simple solution to System} (\ref{v1a}) {\it and let ${\cal S}_{i,t,n} =
\{x_n^{i,t,\sigma}(\lambda): \sigma = 1,\ldots,s_i; \ i = 1,\ldots,n\}$ be the set of jets of
simple solutions $x_n^{i,\sigma}(\lambda)$ ($\sigma = 1,\ldots,s_i; \ i = 1,\ldots,n$) (where $t$
is a number of members of $x_n^{i,t,\sigma}(\lambda)$; $t\ge r_i$; $r_i$ is a defining number of
jets of the solutions $x_n^{i,\sigma}(\lambda)$) for each equation of System (\ref{v2}). Let the
sets ${\cal S}_{i,t,n} = \{x_n^{i,t,\sigma}(\lambda): \sigma = 1,\ldots,s_i; \ i = 1,\ldots,n\}$
have a unique common element $\widetilde{\xi}_n(\lambda)$. Then System} (\ref{nnewsys}) {\it has a
small simple solution $x(\lambda) = (\xi_1(\lambda),\ldots,\xi_{n-1}(\lambda),\xi_n(\lambda))$,
where $\xi_n(\lambda)$ is a simple solution to one of the equations of System (\ref{v2}) whose jet
coincides with $\widetilde{\xi}_n (\lambda)$.}

Applying the Weierstrass preparation theorem and the suitable change of variables we pass from
System (\ref{v1a}) to the consideration of the equivalent system of algebraic equations in the
unknown $x_{n-1}$:
 \begin{equation}\label{newsys1}
 \widetilde{f}^{(n-1)}_j(\lambda,x_1,\ldots,x_{n-1}) = 0 \qquad (j = 1,\ldots,n-1).
 \end{equation}

Let us choose a tree from ${\cal D}_{n-1}$ and use the same scheme for System (\ref{newsys1}).
Following this scheme we get a chain of trees $\delta = (D_n,D_{n-1},\ldots,D_2)$ ($D_j \in {\cal
D}_j$, $j = 2,\ldots,n$) and the corresponding chain of systems (for each system from this chain
the number of equations coincides with the number of unknowns):
 \begin{equation}\label{chain1}
 \begin{array}{c}
 \left\{\begin{array}{l}
 f^{(n)}_1(\lambda,x_1,\ldots,x_{n}) = 0, \\
 \ldots \ldots \ldots \ldots \ldots \ldots \ldots \\
 f^{(n)}_{n}(\lambda,x_1,\ldots,x_{n}) = 0,
 \end{array}\right.
 \stackrel{D_n}{\longrightarrow}
 \left\{\begin{array}{l}
 f^{(n-1)}_1(\lambda,x_1,\ldots,x_{n-1}) = 0, \\
 \ldots \ldots \ldots \ldots \ldots \ldots \ldots \ldots \ldots \\
 f^{(n-1)}_{n-1}(\lambda,x_1,\ldots,x_{n-1}) = 0,
 \end{array}\right.
 \stackrel{D_{n-1}}{\longrightarrow} \ldots \\[24pt] \ldots \stackrel{D_{k}}{\longrightarrow}
 \left\{\begin{array}{l}
 f^{(k-1)}_1(\lambda,x_1,\ldots,x_{k-1}) = 0, \\
 \ldots \ldots \ldots \ldots \ldots \ldots \ldots \ldots \ldots\\
 f^{(k-1)}_{k-1}(\lambda,x_1,\ldots,x_{k-1}) = 0,
 \end{array}\right.
 \stackrel{D_{k-1}}{\longrightarrow} \ldots \stackrel{D_3} {\longrightarrow}
 \left\{\begin{array}{l}
 f^{(2)}_1(\lambda,x_1,x_2) = 0, \\
 f^{(2)}_2(\lambda,x_1,x_2) = 0,
 \end{array}\right.
 \stackrel{D_2}{\longrightarrow} \\[24pt]
 \stackrel{D_2}{\longrightarrow} f_1^{(1)}(\lambda,x_1) = 0.
 \end{array}
 \end{equation}

Let us emphasize that the tree $D_2$ in this chain is defined unequivocally (the set ${\cal D}_2$
consists of one element).

Let us notice that for the system
 \begin{equation}\label{vk}
 \left \{\begin{array}{c}
 f^{(k-1)}_1(\lambda,x_1,\ldots,x_{k-1}) = 0, \\
 \ldots \ldots \ldots \ldots \ldots \ldots \ldots \ldots \ldots \\
 f^{(k-1)}_{k-1}(\lambda,x_1,\ldots,x_{k-1}) = 0.
 \end{array}\right.
 \end{equation}
 it is possible to formulate the statements that are analogues to Lemma 2 and 3.

Let us consider the last system of equations from this chain, i.e. the equation
 \begin{equation}\label{lastscsyst}
 f_1^{(1)}(\lambda,x_1) = 0.
 \end{equation}
Let $x_1(\lambda)$ be a simple solution of this equation. If this solution is a $D_2$-regular
solution then the previous system
 \begin{equation}\label{prelastscsyst}
 \left\{\begin{array}{l}
 f^{(2)}_1(\lambda,x_1,x_2) = 0, \\
 f^{(2)}_2(\lambda,x_1,x_2) = 0,
 \end{array}\right.
 \end{equation}
has a unique simple solution $(x_1(\lambda),x_2(\lambda))$. Proceeding similarly under assumption
of $D_k$-regularities of solution $(x_1(\lambda),\ldots,x_{k-1}(\lambda))$ of the system
 \begin{equation}\label{pk}
 \left\{\begin{array}{l}
 f^{(k-1)}_1(\lambda,x_1,\ldots,x_{k-1}) = 0, \\
 \ldots \ldots \ldots \ldots \ldots \ldots \ldots \ldots \ldots \\
 f^{(k-1)}_{k-1}(\lambda,x_1,\ldots,x_{k-1}) = 0,
 \end{array}\right.
 \end{equation}
we get a simple solution $(x_1(\lambda),\ldots,x_{k-1}(\lambda),x_k(\lambda))$ of the system
 \begin{equation}\label{k}
 \left\{\begin{array}{l}
 f^{(k)}_1(\lambda,x_1,\ldots,x_{k-1},x_{k}) = 0, \\
 \ldots \ldots \ldots \ldots \ldots \ldots \ldots \ldots \ldots \\
 f^{(k)}_{k}(\lambda,x_1,\ldots,x_{k-1},x_{k}) = 0,
 \end{array}\right.
 \end{equation}
where $k = 2,\ldots,n-1$. At last under assumption of $D_n$-regularities of the constructed
solution $(x_1(\lambda),\ldots,x_{n-1}(\lambda))$ of the system
 \begin{equation}\label{pn}
 \left\{\begin{array}{l}
 f^{(n-1)}_1(\lambda,x_1,\ldots,x_{n-1}) = 0, \\
 \ldots \ldots \ldots \ldots \ldots \ldots \ldots \ldots \ldots \\
 f^{(n-1)}_{n-1}(\lambda,x_1,\ldots,x_{n-1}) = 0,
 \end{array}\right.
 \end{equation}
we get a simple solution $(x_1(\lambda),\ldots,x_{n-1}(\lambda),x_n(\lambda))$ of the system
 \begin{equation}\label{n}
 \left\{\begin{array}{l}
 f^{(n)}_1(\lambda,x_1,\ldots,x_{n-1},x_{n}) = 0, \\
 \ldots \ldots \ldots \ldots \ldots \ldots \ldots \ldots \ldots \\
 f^{(n)}_{n}(\lambda,x_1,\ldots,x_{n-1},x_{n}) = 0.
 \end{array}\right.
 \end{equation}

It is obvious that the common number of such different chains is equal to $\prod\limits_{i=2}^n
i^{i-2}$. A solution $x(\lambda) = (x_1(\lambda),\ldots,x_n(\lambda))$ of System (\ref{scsyst}) is
called {\it an effectively computable solution} if $(x_1(\lambda),\ldots,x_{k-1}(\lambda))$ are
$D_k$-regular solutions ($k=1,\ldots,n$).

It is obvious that an effectively computable solution is a simple solution since each component of
this solution is a simple solution of the corresponding system. The scheme of the construction of
effectively computable solutions to System (\ref{scsyst}) will be called the {\it modified
R\"uckert--Lefschetz scheme}.

It is necessary to notice that contrary to the R\"uckert--Lefschetz scheme the modified
R\"uckert--Lefschetz scheme does not allow to get the full description of the solutions to System
(\ref{scsyst}), however in some cases the modified R\"uckert--Lefschetz scheme allows to construct
effectively computable solutions to System (\ref{scsyst}).

The argument presented above implies

{\bf Theorem 4.} {\it The modified R\"uckert--Lefschetz scheme is an effective scheme for
construction of the set of solutions to System (\ref{scsyst}) if and only if this set consists only
of effectively computable solutions.}

{\bf The case of real effectively computable solutions.} Above we supposed that the parameter
$\lambda$ and the unknowns $x_1,\ldots,x_m$ are complex numbers and the coefficients of the
expansion in series $f_j(\lambda,x_1,\ldots,x_m)$ ($j = 1,\ldots,n$) are complex numbers as well.
Therefore effectively computable solutions to System (\ref{scsyst}) constructed by the modified
R\"uckert--Lefschetz scheme in the general case are complex. However, a lot of applications as a
rule represent the cases when the parameter $\lambda$ and the unknowns $x_1,\ldots,x_m$ are real
numbers and the coefficients of expansion $f_j(\lambda,x_1,\ldots,x_m)$ ($j = 1,\ldots,n$) are real
as well. Hereafter we show how to determine in such case which of the effectively computable
solutions to System (\ref{scsyst}) constructed by the modified R\"uckert--Lefschetz scheme are
real.

Let us use Newton's diagram method for the construction of the component $x_k(\lambda)$ ($k =
1,\ldots,n$) of an effectively computable solution $x(\lambda) =
(x_1(\lambda),\ldots,x_n(\lambda))$ to System (\ref{scsyst}).

The Newton's diagram method allows to construct the set of solutions to a scalar equation in the
parameter $\lambda$. Thus each solution to this equation can be presented in the form of series
converging in some neighborhood of zero.

The components $x_k(\lambda)$ ($k = 1,\ldots,n$) of effectively computable solutions $x(\lambda) =
(x_1(\lambda),\ldots,x_n(\lambda))$ to System (\ref{scsyst}) are defined by the scalar equations.
Thus each component $x_k(\lambda)$ ($k = 1,\ldots,n$) is a simple solution of the corresponding
scalar equation and can be represented in the form of converging series in some neighborhood of
zero
 $$\begin{array}{c}
 x_k(\lambda) = \gamma_0\lambda^\frac{\tau_0}{\tau} + \gamma_1\lambda^\frac{\tau_{1}}{\tau} +
 \ldots + \gamma_{l}\lambda^\frac{\tau_{l}}{\tau} + o(\lambda^\frac{\tau_{l}}{\tau}) \\[6pt]
 \quad (\gamma_{i} \in {\Bbb C}, \ \gamma_{0}\neq 0,\ \tau, \tau_{i} \in {\Bbb N}, \ i=0,\ldots,l).
 \end{array}$$
Therefore, if all the members of the expansion $x_k (\lambda)$ were real up to the defining number
$r_k$ then all the subsequent members of expansion $x_k(\lambda)$ will be real as well.

It is necessary to notice that realness of the coefficient $\gamma_{l}$ not means realness of the
member $\gamma_{l}\lambda^\frac{\tau_{l}}{\tau}$ of the expansion $x_k(\lambda)$ ($k=1,\ldots,n$)
since at different values $\lambda$ ($\lambda\ge 0$ and $\lambda < 0$) the conditions on the
realness of the member $\gamma_{l}\lambda^\frac{\tau_{l}}{\tau}$ can be different.

Let us notice also that the members of the expansion of the components
$x_1(\lambda),\ldots,x_n(\lambda)$ of a simple solution $x(\lambda) =
(x_1(\lambda),\ldots,x_n(\lambda))$ to System (\ref{scsyst}) are defined by a finite number of
coefficients in the expansion into series of the left hand parts of the equations of System
(\ref{scsyst}).

The argument prsented above implies

{\bf Theorem 5.} {\it The modified R\"uckert--Lefschetz scheme allows to determine real effectively
compu\-table solutions to System (\ref{scsyst}). Thus an effectively computable solution
$x(\lambda) = (x_1(\lambda),\ldots,x_n(\lambda))$ is real if and only if first $r_k$ members of the
expansion of each component $x_k(\lambda)$ ($k = 1,\ldots,n$) are real.}

\bigskip

\def\refname{{\normalsize \bf Literature}}\label{L}

\end{document}